# Optimal Communication Scheduling and Remote Estimation over an Additive Noise Channel [★]

Xiaobin Gao, Emrah Akyol, Tamer Başar

*University of Illinois at Urbana-Champaign*

**Abstract**

This paper considers a sequential sensor scheduling and remote estimation problem with one sensor and one estimator. The sensor makes sequential observations about the state of an underlying memoryless stochastic process, and makes a decision as to whether or not to send this measurement to the estimator. The sensor and the estimator have the common objective of minimizing expected distortion in the estimation of the state of the process, over a finite time horizon. The sensor is either charged a cost for each transmission or constrained on transmission times. As opposed to the prior work where communication between the sensor and the estimator was assumed to be perfect (noiseless), in this work an additive noise channel with fixed power constraint is considered; hence, the sensor has to encode its message before transmission. Under some technical assumptions, we obtain the optimal encoding and estimation policies in conjunction with the optimal transmission schedule. The impact of the presence of a noisy channel is analyzed numerically based on dynamic programming. This analysis yields some rather surprising results such as a phase-transition phenomenon in the number of used transmission opportunities, which was not encountered in the noiseless communication setting.

## 1 Introduction

### 1.1 Background

The communication scheduling and remote state estimation problem arises in the applications of wireless sensor networks, such as environmental monitoring and networked control systems. As an example of environmental monitoring, researchers at the National Aeronautics and Space Administration (NASA) Earth Science group are interested in monitoring the evolution of the soil moisture, which is used in weather forecast, ecosystem process simulation and the like [1]. In order to achieve that goal, the sensor networks are built over an area of interest. The sensors collect data on the soil moisture and send them to the decision unit at NASA via wireless communication. The decision unit at NASA forms estimates on the evolution of the soil moisture based on the messages received from the sensors. Similarly, in networked control systems, where the objective is to control some remote plants, sensor networks are built to measure the states of the remote plants. Sensors transmit their measurements to the controller via a wireless communication network, and the controller estimates the state of the remote plant and generates a control signal based on that estimate [2]. In both scenarios, the quality of the remote state estimation strongly affects the quality of decision making at the remote site, that is, weather prediction or control signal generation. The networked sensors are usually constrained by limits on power [3]. They are not able to communicate with the estimator at every time

[★] This paper was not presented at any IFAC meeting. This research was supported in part by NSF under grant CCF 11-11342, and in part by the U.S. Air Force Office of Scientific Research (AFOSR) MURI grants FA9550-10-1-0573 and N00014-16-1-2710.

*Email addresses:* `gao16@illinois.edu` (Xiaobin Gao), `akyol@illinois.edu` (Emrah Akyol), `basar1@illinois.edu` (Tamer Başar).

step and thus, the estimator has to produce its best estimate based on the partial information received from the sensors. Therefore, the communication between the sensors and the estimator should be scheduled judiciously, and the estimator should be designed properly, so that the state estimation error is minimized subject to the communication constraints.

*1.2 Literature Review*

Research on the general sensor scheduling problem dates back to the 1970s. In one of the earliest works [4], the problem formulation is such that only one out of several sensors can be selected at each instant of time to observe the output of a linear stochastic system. Using the measurements over a finite time interval, the goal is to form prediction on some future state of the system. Furthermore, each sensor is associated with a certain measurement cost. The author proposed an off-line deterministic sensor scheduling strategy that minimizes the sum of measurement cost over the time interval and prediction error. Gupta *et al* [5] studied the sensor scheduling problem over infinite time horizon. Similar to the problem in [4], only one sensor can be selected at each instant of time. However, there is no measurement cost associated with each sensor. The authors proposed an off-line stochastic sensor scheduling strategy such that the expected steady state estimation error is minimized. Yang and Shi [6] studied the off-line sensor scheduling problem where there is only one sensor observing the state of a linear stochastic system. The sensor can communicate with the remote estimator only a limited number of times. The objective was to minimize the cumulative estimation error over a finite time horizon. It was shown that the optimal sensor scheduling strategy is to distribute the limited communication opportunities uniformly over the time horizon. The authors of the papers discussed above considered off-line sensor scheduling problems. "Off-line sensor scheduling" means the sensor is scheduled to take observation or conduct communication based on some a priori information about the system (e.g. statistics of random variables, system matrices). The on-line information (e.g. sensor's observation, battery's energy level) is not taken into account when making schedules. Some other selected work on off-line sensor scheduling problems can be found in [7–9].

With the advances in hardware devices, sensors are endowed with stronger computational capabilities. Consequently, the sensors are able to make schedules based on all the information they have (a priori information as well as on-line information), which motivates the formulation of on-line sensor scheduling problems. Åström and Bernhardsson [10] considered a state estimation problem with a first-order stochastic system. They compared the estimation error over infinite time horizon obtained by periodic sampling and threshold event-triggered sampling. The periodic sampling is one of the off-line sensor scheduling strategies while the threshold event-triggered sampling is one of the on-line sensor scheduling strategies. They showed that the threshold event-triggered sampling, which is also called "threshold-based communication strategy", leads to better performance in state estimation compared with periodic sampling. The global optimality of threshold-based communication strategy in this contexts is proved later by Nar and Başar [11]. Imer and Başar [12] considered the on-line sensor scheduling and remote state estimation problem over a finite time horizon. In the formulation, the sensor is restricted to communicate only a limited number of times. By considering the communication strategies within the class of threshold-based strategies, the paper has shown that there exists a unique threshold-based communication strategy achieving the best performance on remote state estimation. Furthermore, the optimal threshold can be computed by solving a dynamic programming equation. Bommannavar and Başar [13] later extended the result of [12] to multi-dimensional systems. The continuous-time version of the problem in [12] has been studied by Rabi *et al* [14]. Xu and Hespanha [15] considered the networked control problem involving state estimation and communication scheduling, which can be viewed as a sensor scheduling and remote estimation problem. They fixed the estimator to be Kalman-like and designed an event-triggered sensor that minimizes the time average of the sum of the communication cost and estimation error over infinite time horizon. They showed that the optimal communication strategy is deterministic and stationary, and is a function of the estimation error. Wu *et al* [16] considered the sensor scheduling and estimation problem subject to constraints on the average communication rate over infinite time horizon. The authors assumed that the sensor has noisy observations on the system state. By restricting the sensor scheduling strategies to the threshold event-triggered class, they derived the exact minimum mean square error (MMSE) estimator. However, the exact MMSE estimator is nonlinear and thus computationally intractable. Under a Gaussian assumption on the a priori distribution, the authors derived an approximate MMSE estimator, which is Kalman-like. Based on the approximated MMSE estimator, the authors derived conditions on the thresholds so that the average sensor communication rate will not exceed its upper bound. You and Xie [17] extended the work in [16] by deriving conditions on the thresholds so that the estimator is stable. Han *et al.* [18] showed that if the sensor is fixed to apply some stochastic event-triggered strategy, then the exact MMSE estimator is Kalman-like. Other selected work on remote estimation with event-based sensor operations can be found in [19, 20]. The work in [16–18] can also be viewed as Kalman-filtering with scheduled observations, which is related to Kalman-filtering with intermittent observations studied in [21, 22].

The approaches of [15, 16] involved fixing the communication strategies or estimation strategies to be of a certain type and then deriving the corresponding optimal estimation strategies and communication strategies, respectively. The approach of [12], on the other hand, is to derive the jointly optimal communication strategies and estimation strategies. Similarly, Lipsa and Martins [23] considered the sensor scheduling and remote estimation problem where the sensor is not constrained by communication times but is charged a communication cost. They proposed a threshold event-triggered sensor and a Kalman-like estimator and proved that the proposed sensor and estimator are jointly optimal, minimizing the sum of communication cost and estimation error over a finite time horizon. Nayyar *et al* [24] considered a similar problem where the sensor is equipped with an energy harvesting sensor. In the work of [24], the problem formulation is such that the sensor is constrained by the energy level of the battery and is also charged a communication cost. It is shown in [24] that an energy dependent threshold event-triggered sensor and a Kalman-like estimator are jointly optimal. Hence, the result of [24] can be viewed as generalization of the results of [12, 23]. In both [23] and [24], majorization theory was used to prove the optimality of the respective results, which is closely related to the approach in [25].

It is worth drawing attention to the two different types of constraints that arise in the works mentioned above–hard and soft constraints–as featured in the problem setups of [12] and [23]. In the problem of [12], the sensor can only communicate for a pre-specified number of times. Such a communication constraint is called *hard constraint*. In the work of [23], however, the sensor is charged a communication cost. This kind of communication constraint is called *soft constraint*. In the problem with hard constraint, the communication strategy must take the remaining communication opportunities as a variable and schedule no communication if there is no remaining opportunity. Such communication strategies guarantee that the number of transmissions made over the time horizon of interest will not exceed the given constraint. In the problem with soft constraint, however, the sensor is not constrained by the number of transmissions, therefore the communication strategy need not take the remaining communication opportunities (which is always equal to the remaining time steps) as a variable. Therefore, the results obtained in one problem cannot be applied directly to the other. For example, if we apply the communication strategy obtained in [23] to the problem of [12], then there exists a positive probability that the sensor decides to communicate at every instant of time, which certainly would violate the hard constraint on communication times. A detailed discussion on the difference between optimization problems with soft constraints and hard constants can be found in [26, 27].

*1.3 Contributions*

In this paper, we extend the lines of research in [12], [23], and [24]. In previous works, the communication between the sensor and the estimator was assumed to be perfect (no channel noise), which may not be realistic, even though it was an important first step. This paper investigates the effect of communication channel noise on the design of optimal sensor scheduling and remote estimation strategies. The paper consists of two parts: In the first part, we consider a discrete time sensor scheduling and remote estimation problem with a soft constraint. At each time step, the sensor makes a perfect observation on the state of an independent identically distributed (i.i.d.) source. Next, the sensor decides whether to transmit its observation to the remote estimator or not. The sensor has a soft communication constraint (i.e., the sensor is charged a cost for each transmission). Since the communication channel is noisy, the sensor encodes the message before transmitting it to the estimator. The remote estimator generates a real-time estimate on the state of the source based on the noise corrupted messages received from the sensor. The estimator is charged for estimation error. Our goal is to design the communication scheduling strategy and encoding strategy for the sensor, and the estimation strategy (decoding strategy) for the estimator, to minimize the expected value of the sum of communication cost and estimation cost over the time horizon. Our solution consists of a threshold-based communication scheduling strategy, and a pair of piecewise linear encoding/decoding strategies. We prove optimality under some technical assumptions. Then, we extend the result to the problem with a hard constraint. Using a dynamic programming approach, we obtain the optimal communication scheduling, encoding and decoding strategies. Beyond the qualititively expected results, we notice some rather surprising effects of the noisy communication considerations in this class of remote estimation problems. For example, over a time horizon $T$ and with a hard transmission limit, $N \leq T$, if the state realizations were so that at time step $K$, the sensor has used only $N - T + K$ transmissions out of $N$, the intuitively appealing solution to the noiseless variation of the problem was to transmit from that time on all the observed state realizations without any thresholding, i.e., the threshold is effectively set to zero for samples at time steps $K+1, \ldots, T$. However, in the noisy setting, we have noticed that this is not the case, the sensor may not use all the transmission opportunities left. This is due to the fact that threshold information–that is whether or not the state sample belongs to an interval– may be more valuable than a "noisy" observation of the state. In fact, depending on the signal-to-noise ratio (SNR) of the channel, there is a fixed number of useful (in average) number of transmissions, and allowing transmissions more than this number, on the average, does not help decrease the expected mean square error (MSE).

*1.4 Organization*

The rest of the paper is organized as follows: in Section 2, we formulate the optimization problems with soft/hard constraints. In Sections 3 and 4, we present the main results for problems with soft/hard constraints. In Section 5, we present some numerical results for the problem with hard constraint. Finally, in Section 6, we draw concluding remarks and discuss future work.

The authors have five conference papers on the general topic of this paper, listed as [31]-[35]. The specific topics and results of the last three, that is [33]-[35], are beyond the scope of this paper, as explained in Section VI (Conclusions). The first two, that is [31] and [32], have some overlap as far as the problem formulations go, but the current paper substantially improves upon the results in these two papers, as explained in Remark 9 in Section IV.

## 2 Problem Formulation

*2.1 System Model*

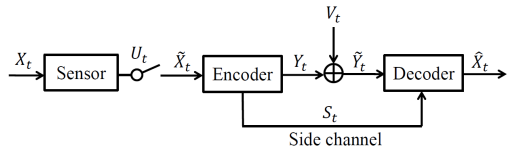

Fig. 1. System model

Consider a discrete time communication scheduling and remote estimation problem over a finite-time horizon, that is, $t = 1, 2, \ldots, T$. There is one sensor, one encoder and one remote estimator (which is also called "decoder"). A source process $\{X_t\}$ is a one-dimensional, independent, and identically distributed (i.i.d.) stochastic process, which has density $p_X$. At time $t$, the sensor observes $X_t$. Since the sensor is assumed to have communication constraint (which will be introduced later), it needs to decide whether or not to transmit its observation. Let $U_t \in \{0, 1\}$ be the sensor's decision at time $t$, where $U_t = 1$ stands for transmission and $U_t = 0$ stands for no transmission. The communication channel is assumed to be noisy. Hence, if the sensor decides to transmit its observation, it sends $X_t$ to the encoder. If the sensor decides not to transmit, it does not send anything to the encoder but a free symbol $\epsilon$ stands for its decision. Denote by $\tilde{X}_t$ the message received by the encoder; then

$$\tilde{X}_t = \begin{cases} X_t, & \text{if } U_t = 1 \\ \epsilon, & \text{if } U_t = 0 \end{cases}$$

If the encoder receives $X_t$ from the sensor, it sends an encoded message $Y_t$ to the communication channel. The encoder operates under the average power constraint:

$$\mathbb{E}[Y_t^2 | U_t = 1] \leq P_T,$$

where the expectation is taken over $Y_t$. Furthermore, $P_T$ is known and is invariant of time. The encoded message $Y_t$ is corrupted by an additive channel noise $V_t$. The noise process $\{V_t\}$ is a one-dimensional i.i.d. stochastic process with density $p_V$. When sending $Y_t$ to the communication channel, the encoder is able to transmit the sign of $X_t$ to the decoder via a side channel, which is assumed to be noise-free. If the encoder receives $\epsilon$ from the sensor, it does not send anything to the communication channel, while it sends $\epsilon$ to the side channel. Consequently, the decoder can deduce the sensor's decision from the message conveyed via the side-channel. We use $\tilde{Y}_t$ and $S_t$ to denote the messages received by the decoder from the communication channel and the side channel, respectively, that is

$$\tilde{Y}_t = \begin{cases} Y_t + V_t, & \text{if } U_t = 1 \\ V_t, & \text{if } U_t = 0 \end{cases}, \quad S_t = \begin{cases} \text{sgn}(X_t), & \text{if } U_t = 1 \\ \epsilon, & \text{if } U_t = 0 \end{cases}$$

After receiving $\tilde{Y}_t$ and $S_t$, the decoder produces an estimate on $X_t$, denoted by $\hat{X}_t$. The estimator will be charged for distortion in estimation. We assume that the distortion function $\rho(X_t, \hat{X}_t)$ is the squared error, $(X_t - \hat{X}_t)^2$.

## 2.2 Communication Constraint

The sensor is said to have a *soft constraint* if there is a positive cost function associated with $U_t$, denoted by $\mathcal{C}(U_t)$. Here, the cost function is assumed to have the form of

$$\mathcal{C}(U_t) = cU_t = \begin{cases} 0, & \text{if } U_t = 0 \\ c, & \text{if } U_t = 1 \end{cases},$$

where $c$ is called the communication cost ($c > 0$), which is known and is invariant of time. The sensor is said to have a *hard constraint* if it is restricted to use the noisy channel for no more than $N$ times, i.e.,

$$\sum_{t=1}^{T} \mathbf{1}_{\{U_t = 1\}} \leq N,$$

where $\mathbf{1}_{\{\cdot\}}$ is the indicator function, and $0 < N < T$.

## 2.3 Decision Strategies

Assume that at time $t$, the sensor has memory on all its observations up to $t$, denoted by $X_{1:t}$, and all the decisions it has made up to $t-1$, denoted by $U_{1:t-1}$. The sensor determines whether or not to transmit its observation at time $t$, based on its current information $(X_{1:t}, U_{1:t-1})$, namely

$$U_t = f_t(X_{1:t}, U_{1:t-1}),$$

where $f_t$ is the communication scheduling policy at time $t$, and $\mathbf{f} = \{f_1, f_2, \ldots, f_T\}$ is the communication scheduling strategy.

Similarly, at time $t$, the encoder has memory on all the messages received from the sensor up to $t$, denoted by $\tilde{X}_{1:t}$, and all the messages it has sent to the communication channel and the side channel up to $t-1$, denoted by $Y_{1:t-1}$ and $S_{1:t-1}$, respectively. The encoder generates the encoded message at time $t$, based on its current information $(\tilde{X}_{1:t}, Y_{1:t-1}, S_{1:t-1})$, namely

$$Y_t = g_t(\tilde{X}_{1:t}, Y_{1:t-1}, S_{1:t-1}),$$

where $g_t$ is the encoding policy at time $t$, and $\mathbf{g} = \{g_1, g_2, \ldots, g_T\}$ is the encoding strategy.

Finally, we assume that at time $t$, the decoder has memory on all the messages received from the communication channel up to $t$, denoted by $\tilde{Y}_{1:t}$, and all the messages received from the side channel up to $t$, which are $S_{1:t}$. The decoder generates the estimate at time $t$, based on its current information $(\tilde{Y}_{1:t}, S_{1:t})$, namely

$$\hat{X}_t = h_t(\tilde{Y}_{1:t}, S_{1:t}),$$

where $h_t$ is the decoding policy at time $t$, and $\mathbf{h} = \{h_1, h_2, \ldots, h_T\}$ is the decoding strategy.
**Remark 1.** *Although we do not assume that the decoder has memory on its previous estimates up to $t$, yet it can deduce them from $(\tilde{Y}_{1:t-1}, S_{1:t-1})$ and $h_1, h_2, \ldots, h_{t-1}$.*

For simplicity, we call the sensor, the encoder, and the decoder as "decision makers". Correspondingly, we call the communication scheduling policy, encoding policy, and decoding policy as "decision policies".

## 2.4 Optimization Problem

Consider the settings described above, with the time horizon $T$, the probability density functions $p_X$ and $p_V$, and the power constraint $P_T$ as given.

**Optimization problem with soft constraint**: Given the communication cost $c$, determine $(\mathbf{f}, \mathbf{g}, \mathbf{h})$ minimizing the function

$$J(\mathbf{f}, \mathbf{g}, \mathbf{h}) = \mathbb{E}\left\{\sum_{t=1}^{T} cU_t + (X_t - \hat{X}_t)^2\right\}.$$

**Optimization problem with hard constraint**: Given the number of transmission opportunities $N < T$, determine $(\mathbf{f}, \mathbf{g}, \mathbf{h})$ minimizing, under the hard constraints, the function

$$J(\mathbf{f}, \mathbf{g}, \mathbf{h}) = \mathbb{E}\left\{\sum_{t=1}^{T} (X_t - \hat{X}_t)^2\right\}.$$

## 3  Optimization Problem with Soft Constraint

To begin with, we show that the optimization problem with soft constraint can be simplified to a "one-stage" problem, as described in the theorem below.

**Theorem 1.** *Consider the optimization problem in Section 2-D with the soft constraint.*

(1) *Without loss of optimality, one can restrict all the decision makers to apply the decision policies $(f_t, g_t, h_t)$ in the forms of:*
$$U_t = f_t(X_t),\ Y_t = g_t(\tilde{X}_t),\ \hat{X}_t = h_t(\tilde{Y}_t, S_t). \tag{1}$$

(2) *Without loss of optimality, one can restrict all the decision makers to apply stationary decision strategies $(\boldsymbol{f}, \boldsymbol{g}, \boldsymbol{h})$, i.e.,*
$$\boldsymbol{f} = \{f, f, \ldots, f\},\quad \boldsymbol{g} = \{g, g, \ldots, g\},\quad \boldsymbol{h} = \{h, h, \ldots, h\}.$$

Before proving Theorem 1, we first introduce the following notations. For any $1 \leq a \leq b \leq T$, let $\mathbf{f}_{a:b}, \mathbf{g}_{a:b}, \mathbf{h}_{a:b}$ denote the subsets of $\mathbf{f}, \mathbf{g}, \mathbf{h}$ such that

$$\mathbf{f}_{a:b} = \{f_a, f_{a+1}, \ldots, f_b\},\quad \mathbf{g}_{a:b} = \{g_a, g_{a+1}, \ldots, g_b\},\quad \mathbf{h}_{a:b} = \{h_a, h_{a+1}, \ldots, h_b\}.$$

Furthermore, let $I_{st}, I_{et}, I_{dt}$ denote the information about the past system states available to the sensor, the encoder, and the decoder, respectively, at time $t$ $(t > 1)$, i.e.,

$$I_{st} = \{X_{1:t-1}, U_{1:t-1}\},\quad I_{et} = \{\tilde{X}_{1:t-1}, Y_{1:t-1}, S_{1:t-1}\},\quad I_{dt} = \{\tilde{Y}_{1:t-1}, S_{1:t-1}\}.$$

Furthermore, let $I_t$ be union of $I_{st}$, $I_{et}$, and $I_{dt}$, i.e.,

$$I_t = \{X_{1:t-1}, U_{1:t-1}, \tilde{X}_{1:t-1}, Y_{1:t-1}, S_{1:t-1}, \tilde{Y}_{1:t-1}\}.$$

*Proof of Theorem 1*: It is easy to see the validity of the following sequence of equalities:

$$\min_{\mathbf{f}, \mathbf{g}, \mathbf{h}} J(\mathbf{f}, \mathbf{g}, \mathbf{h})$$
$$= \min_{\mathbf{f}, \mathbf{g}, \mathbf{h}} \mathbb{E}\left\{\sum_{t=1}^{T} cU_t + (X_t - \hat{X}_t)^2\right\}$$
$$= \min_{f_1, g_1, h_1} \mathbb{E}\left\{cU_1 + (X_1 - \hat{X}_1)^2 + \min_{\mathbf{f}_{2:T}, \mathbf{g}_{2:T}, \mathbf{h}_{2:T}} \mathbb{E}\left\{\sum_{t=2}^{T} cU_t + (X_t - \hat{X}_t)^2\right\}\right\}$$
$$= \min_{f_1, g_1, h_1} \mathbb{E}\left\{cU_1 + (X_1 - \hat{X}_1)^2 + \min_{f_2, g_2, h_2} \mathbb{E}\left\{cU_2 + (X_2 - \hat{X}_2)^2 + \ldots + \min_{f_T, g_T, h_T} \mathbb{E}\left\{cU_T + (X_T - \hat{X}_T)^2\right\}\ldots\right\}\right\}$$

Then, at time $t = T$, the optimization problem is to design $(f_T, g_T, h_T)$ minimizing

$$J_{T_1}(f_T, g_T, h_T) := \mathbb{E}\left\{cU_T + (X_T - \hat{X}_T)^2\right\},$$

call it *Problem T1*. Recall that the decisions at time $T$ are generated by

$$U_T = f_T(X_T, I_{sT}), \quad Y_T = g_T(\tilde{X}_T, I_{eT}), \quad \hat{X}_T = h_T(\tilde{Y}_T, S_T, I_{dT}).$$

We will now show that using information about the past $(I_{sT}, I_{eT}, I_{dT})$ when making decisions cannot help improve the performance (that is, reduce the expected cost). Consider another problem, call it *Problem T2*, where $I_T$ is available to all the decision makers, and one needs to design $(f'_T, g'_T, h'_T)$ minimizing

$$J_{T_2}(f'_T, g'_T, h'_T) := \mathbb{E}\left\{cU_T + (X_T - \hat{X}_T)^2\right\},$$

where

$$U_T = f'_T(X_T, I_T), \quad Y_T = g'_T(\tilde{X}_T, I_T), \quad \hat{X}_T = h'_T(\tilde{Y}_T, S_T, I_T).$$

Since the sensor, the encoder, and the decoder can always ignore the redundant information and behave as if they only know $I_{sT}, I_{eT}, I_{dT}$, respectively, the performance of the system in *Problem T2* is no worse than that in *Problem T1*, i.e.,

$$\min_{(f'_T, g'_T, h'_T)} J_{T_2}(f'_T, g'_T, h'_T) \leq \min_{(f_T, g_T, h_T)} J_{T_1}(f_T, g_T, h_T).$$

Similarly, consider a third problem, call it *Problem T3*, where $I_{sT}, I_{eT}, I_{dT}$ are not available to the sensor, the encoder, and the decoder, respectively. One needs to design $(f''_T, g''_T, h''_T)$ to minimize

$$J_{T_3}(f''_T, g''_T, h''_T) = \mathbb{E}\left\{cU_T + (X_T - \hat{X}_T)^2\right\},$$

where

$$U_T = f''_T(X_T), \quad Y_T = g''_T(\tilde{X}_T), \quad \hat{X}_T = h''_T(\tilde{Y}_T, S_T)$$

By a similar argument as above, the system in *Problem T1* cannot perform worse than the system in *Problem T3*. Hence,

$$\min_{(f_T, g_T, h_T)} J_{T_1}(f_T, g_T, h_T) \leq \min_{(f''_T, g''_T, h''_T)} J_{T_3}(f''_T, g''_T, h''_T).$$

Let us come back to *Problem T2*. One can observe that the communication cost $c$, the distortion function $\rho(\cdot, \cdot)$, and the power constraint of the encoder do not depend on $I_T$. Furthermore, since $\{X_t\}$ and $\{V_t\}$ are i.i.d. stochastic processes, $X_T$ and $V_T$ are also independent of $I_T$. Therefore, there is no loss of optimality in restricting

$$U_T = f'_T(X_T), Y_T = g'_T(\tilde{X}_T), \hat{X}_T = h'_T(\tilde{Y}_T, S_T),$$

and thus

$$\min_{(f'_T, g'_T, h'_T)} J_{T_2}(f'_T, g'_T, h'_T) = \min_{(f''_T, g''_T, h''_T)} J_{T_3}(f''_T, g''_T, h''_T)$$

The equality above indicates that in *Problem T1* the sensor, the encoder, and the decoder can ignore their information about the past, namely $I_{sT}, I_{eT}$, and $I_{dT}$, respectively, and there is no loss of optimality in restricting

$$U_T = f_T(X_T), Y_T = g_T(\tilde{X}_T), \hat{X}_T = h_T(\tilde{Y}_T, S_T).$$

Since $(f_T, g_T, h_T)$ do not take $I_T$ as a parameter, the design of $(f_T, g_T, h_T)$ is independent of the design of $(f_{1:T-1}, g_{1:T-1}, h_{1:T-1})$. Consequently, the problem can be viewed as a $T-1$ stage problem and a one stage problem. By induction, we can show that $(f_1, g_1, h_1), (f_2, g_2, h_2), \ldots, (f_T, g_T, h_T)$ can be designed independently, and $(f_t, g_t, h_t)$ is designed to minimize the stage-wise cost $\mathbb{E}\{cU_t + (X_t - \hat{X}_t)^2\}$. Hence, the optimal decision policies $(f_t, g_t, h_t)$ are in the form of (1). Furthermore, since $\{X_t\}$ and $\{V_t\}$ are i.i.d. stochastic processes, the optimal decision policies $(f_t, g_t, h_t)$ should be the same for all $t = 1, 2, \ldots, T$. □

By Theorem 1, the optimization problem with soft constraint can be reduced to a "one-stage" problem. Therefore, for simplicity we suppress the subscript for time in all the expressions for the rest of this section. To present our main results for the one-stage problem, we need the following four assumptions:

**Assumption 1.** *The source density $p_X$ is nonatomic, even, and log-concave with support $\mathbb{R}$. Furthermore, $p_X$ is continuously differentiable on $(0, \infty)$ (and on $(-\infty, 0)$ by symmetry).*

**Remark 2.** *There are several probability density functions satisfying Assumption 1, e.g., Gaussian distribution with zero mean, Laplace distribution and a few others. For simplicity, we assume here that $p_X$ has support $\mathbb{R}$. However, the results also hold for the source density with support $(-a, a), a > 0$, e.g., uniform distribution. In this case, we require that $p_X$ is continuously differentiable on $(0, a)$.*

Given any communication scheduling policy $f$, let $\mathcal{T}_0^f$, $\mathcal{T}_{1+}^f$, and $\mathcal{T}_{1-}^f$ be the *non-transmission region*, the *positive transmission region* and the *negative transmission region*, respectively, where

$$\mathcal{T}_0^f := \{x \in \mathbb{R} | f(x) = 0\}, \ \ \mathcal{T}_{1+}^f := \{x > 0 | f(x) = 1\}, \ \ \mathcal{T}_{1-}^f := \{x < 0 | f(x) = 1\}$$

Note that $\mathcal{T}_0^f, \mathcal{T}_1^f, \mathcal{T}_2^f$ may not be connected regions. Then, we make the following assumption on the communication scheduling policy.

**Assumption 2.** *The sensor is restricted to apply the communication scheduling policy $f$ satisfying*

$$\mathbb{E}[X | X \in \mathcal{T}_{1-}^f] < \mathbb{E}[X | X \in \mathcal{T}_0^f] < \mathbb{E}[X | X \in \mathcal{T}_{1+}^f] \tag{2}$$

**Remark 3.** *There is a wide class of communication scheduling policies satisfying inequality (2). For example, given any even communication scheduling policy $f$, i.e.,*

$$f(x) = f(-x) \in \{0, 1\},$$

*and any even source density function $p_X$, we have*

$$\mathbb{E}[X | X \in \mathcal{T}_{1-}^f] < 0, \ \ \mathbb{E}[X | X \in \mathcal{T}_0^f] = 0, \ \ \mathbb{E}[X | X \in \mathcal{T}_{1+}^f] > 0,$$

*Then, Assumption 2 is satisfied.*

**Assumption 3.** *The communication channel noise $V$ has zero mean, and finite variance, denoted by $\sigma_V^2$. Furthermore, $V$ and $X$ are independent conditioned on $\mathrm{sgn}(X)$.*

**Assumption 4.** *The encoder and the decoder are restricted to apply piecewise affine policies:*

$$g(\tilde{X}) = \begin{cases} S \cdot \alpha(S) \cdot (X - \mathbb{E}[X | U = 1, S]), & \text{if } U = 1 \\ 0, & \text{if } U = 0 \end{cases}$$

$$h(\tilde{Y}, S) = \begin{cases} S \cdot \frac{1}{\alpha(S)} \frac{\gamma}{\gamma + 1} \tilde{Y} + \mathbb{E}[X | U = 1, S], & \text{if } U = 1 \\ \mathbb{E}[X | U = 0], & \text{if } U = 0 \end{cases}$$

*where $\gamma = P_T / \sigma_V^2$, $\alpha(S) = \sqrt{P_T / \mathrm{Var}(X | U = 1, S)}$, and $\mathrm{Var}(X | U = 1, S)$ is the conditional variance.*

One can see that, when $U = 1$, the encoder can deduce $(X, S)$ from $\tilde{X}$, namely $X = \tilde{X}$ and $S = \mathrm{sgn}(\tilde{X})$. Furthermore, when applying the encoding policy described above, the power consumption of the encoder can be computed as

$$\begin{aligned}
\mathbb{E}[g^2(\tilde{X}) | U = 1] &= \sum_{i \in \{-1, 1\}} \mathbb{E}[g^2(\tilde{X}) | U = 1, S = i] \, \mathbb{P}(S = i) \\
&= \sum_{i \in \{-1, 1\}} \mathbb{E}\left[\alpha^2(S)(X - \mathbb{E}[X | U = 1, S = i])^2 | U = 1, S = i\right] \mathbb{P}(S = i) \\
&= \sum_{i \in \{-1, 1\}} \frac{P_T}{\mathrm{Var}(X | U = 1, S = i)} \mathbb{E}\left[(X - \mathbb{E}[X | U = 1, S = i])^2 | U = 1, S = i\right] \mathbb{P}(S = i) \\
&= P_T,
\end{aligned}$$

which satisfies the average power constraint. Moreover, the events $U = 0$, $(U = 1, S = -1)$, and $(U = 1, S = 1)$ are equivalent to the events $X \in \mathcal{T}_0^f$, $X \in \mathcal{T}_-^f$, and $X \in \mathcal{T}_{1+}^f$, respectively. Therefore, the encoding and decoding policies $(g, h)$ are determined by the source density $p_X$ and the communication scheduling policy $f$. For simplicity, we use $J(f)$ instead of $J(f, g, h)$ to denote the cost function in the rest of this section.

**Remark 4.** *Note that the assumption of piece-wise affine encoding policies originates from a prior work [28], which analyzed a memoryless zero-sum jamming game between a pair of transmitter and receiver and an adversary that generates an additive channel noise subject to second order (power) statistical constraints. It was shown in [28] that the saddle-point equilibrium associated with this zero-sum game is achieved by linear encoding/decoding policies for the transmitter & receiver pair. Here, we utilize such piece-wise linear policies, not only because they facilitate a tractable analysis but also because they posess such mini-max robustness properties (see [28] for more details).*

**Remark 5.** *Since the communication channel noise $V$ depends on the sign of $X$, the encoder and the decoder will apply different encoding/decoding policies for positive/negative realizations of the source.*

**Theorem 2.** *Consider the one-stage problem under Assumptions 1-4. Then, the optimal communication scheduling policy is of the symmetric threshold type:*

$$f(x) = \begin{cases} 0, & \text{if } |x| < \beta \\ 1, & \text{otherwise} \end{cases},$$

*where $\beta > 0$ is the threshold. Furthermore, there exists a unique value $\beta^*$ minimizing the cost function $J(f)$ among all such thresholds.*

To prove Theorem 2, we need the following definitions and lemmas. We first introduce a quantization problem.

**Quantization Problem**: The problem is one of quantizing the realizations (denoted by $x$) of a real-valued random variable (denoted by $X$) to $N$ codepoints ($N$ is finite and is known) according to some quantization rule (or quantizer) $Q$, i.e,

$$Q(x) = q_i, \text{ if } x \in S_i, \quad i \in \{1, 2, \ldots, N\},$$

where $S_1, S_2, \ldots, S_N$ are called quantization regions and $q_1, q_2, \ldots, q_N$ are the corresponding codepoints. Note that $S_1, S_2, \ldots, S_N$ are mutually disjoint sets and their union equals $\mathbb{R}$. The distortion error between a realization $x$ and the its quantized value $Q(x)$ is $\rho(|x - Q(x)|)$, where $\rho : [0, \infty) \to [0, \infty)$ is called the distortion function. The performance of the quantizer $Q$ is evaluated by its mean distortion error, denoted by $D(Q)$, i.e.,

$$D(Q) := \mathbb{E}\Big[\rho\big(|X - Q(X)|\big)\Big]$$

Then, given the probability distribution of $X$, the optimization problem is to design a quantizer $Q = Q^*$ (i.e., design $\{S_1, S_2, \ldots, S_N\}$ and $\{q_1, q_2, \ldots, q_N\}$) that minimizes $D(Q)$.

We recall here a result on the regularity of the optimal quantizer, which we will use shortly.

**Lemma 1** ([29], Theorem 1 and Corollary 1). *Assume that the source $X$ has nonatomic distribution $p_X$, and $\rho : [0, \infty) \to [0, \infty)$ is convex and nondecreasing. Then, for any $N$-level quantizer $Q$ with quantization regions $\{S_1, S_2, \ldots, S_N\}$ and the corresponding codepoints $\{q_1, q_2, \ldots, q_N\}$, there exists a quantizer $\hat{Q}$ with quantization regions $\{\hat{S}_1, \hat{S}_2, \ldots, \hat{S}_N\}$ and the corresponding codepoints $\{\hat{q}_1, \hat{q}_2, \ldots, \hat{q}_N\}$ such that*

*(1) $\hat{S}_i$ is convex, and $\mathbb{P}(X \in \hat{S}_i) = \mathbb{P}(X \in S_i)$, for all $i = 1, \ldots, N$.*
*(2) If $q_i < q_j$, then $\hat{S}_i < \hat{S}_j$, i.e., $x < y$ for any $x \in \hat{S}_i$ and $y \in \hat{S}_j$.*
*(3) $\hat{q}_i \in \hat{S}_i$, for all $i = 1, \ldots, N$.*
*(4) $D(\hat{Q}) \leq D(Q)$.*

Now returning to our problem, for any communication scheduling policy $f$, we can construct a three level quantizer, denoted by $Q^f$, with quantizing regions $(\mathcal{T}_0^f, \mathcal{T}_{1+}^f, \mathcal{T}_{1-}^f)$ and the corresponding codepoints $(\mathbb{E}[X|X \in \mathcal{T}_0^f], \mathbb{E}[X|X \in \mathcal{T}_{1+}^f], \mathbb{E}[X|X \in \mathcal{T}_{1-}^f])$. Let $D(Q^f)$ be the mean squared distortion of $Q^f$, i.e.,

$$D(Q^f) = \mathbb{E}\Big[(X - Q^f(X))^2\Big] = \sum_{i \in \{0, 1+, 1-\}} \mathbb{E}\Big[(X - \mathbb{E}[X|X \in \mathcal{T}_i^f])^2 | X \in \mathcal{T}_i^f\Big] \mathbb{P}(X \in \mathcal{T}_i^f)$$

$$= \sum_{i \in \{0, 1+, 1-\}} \text{Var}(X|X \in \mathcal{T}_i^f) \mathbb{P}(X \in \mathcal{T}_i^f)$$

By Lemma 1, we have the following result.

**Lemma 2.** *Suppose the source density $p_X$ is nonatomic and even. Then, for any communication scheduling policy $f$ satisfying Assumption 2, we can construct a threshold-based communication scheduling policy $f^{(1)}$ such that*

(1) $\mathcal{T}_0^{f^{(1)}} = (-\beta_2, \beta_1), \mathcal{T}_{1+}^{f^{(1)}} = (\beta_1, \infty), \mathcal{T}_{1-}^{f^{(1)}} = (-\infty, -\beta_2)$, where $\beta_1 > 0, \beta_2 > 0$ are thresholds.

(2) $\mathbb{P}(X \in \mathcal{T}_i^{f^{(1)}}) = \mathbb{P}(X \in \mathcal{T}_i^f)$, for all $i \in \{0, 1+, 1-\}$.

(3) $D(Q^{f^{(1)}}) \leq D(Q^f)$.

*Proof.* By Lemma 1, given a 3-level quantizer $Q^f$ with quantization regions $(\mathcal{T}_0^f, \mathcal{T}_{1+}^f, \mathcal{T}_{1-}^f)$ and the corresponding codepoints $(\mathbb{E}[X|X \in \mathcal{T}_0^f], \mathbb{E}[X|X \in \mathcal{T}_{1+}^f], \mathbb{E}[X|X \in \mathcal{T}_{1-}^f])$, there exists a 3-level quantizer $\hat{Q}$ with quantization regions $(\mathcal{T}_0^{f^{(1)}}, \mathcal{T}_{1+}^{f^{(1)}}, \mathcal{T}_{1-}^{f^{(1)}})$ and corresponding codepoints $(\hat{c}_0, \hat{c}_{1+}, \hat{c}_{1-})$ such that $D(\hat{Q}) \leq D(Q^f)$. Furthermore, $\mathcal{T}_i^{f^{(1)}}$ is convex and $\mathbb{P}(X \in \mathcal{T}_i^{f^{(1)}}) = \mathbb{P}(X \in \mathcal{T}_i^f)$, for all $i \in \{0, 1+, 1-\}$. Since $\mathbb{E}[X|X \in \mathcal{T}_{1-}^f] < \mathbb{E}[X|X \in \mathcal{T}_0^f] < \mathbb{E}[X|X \in \mathcal{T}_{1+}^f]$, we have $\mathcal{T}_{1-}^{f^{(1)}} < \mathcal{T}_0^{f^{(1)}} < \mathcal{T}_{1+}^{f^{(1)}}$. Note that since $\mathcal{T}_{1+}^f \subseteq (0, \infty)$, $\mathcal{T}_{1-}^f \subseteq (-\infty, 0)$, and the source density $p_X$ is even, we have

$$\mathbb{P}(X \in \mathcal{T}_{1+}^{f^{(1)}}) = \mathbb{P}(X \in \mathcal{T}_{1+}^f) \leq \frac{1}{2}, \quad \mathbb{P}(X \in \mathcal{T}_{1-}^{f^{(1)}}) = \mathbb{P}(X \in \mathcal{T}_{1-}^f) \leq \frac{1}{2},$$

which implies that $\mathcal{T}_0^{f^{(1)}} = (-\beta_2, \beta_1), \mathcal{T}_{1+}^{f^{(1)}} = (\beta_1, \infty)$, and $\mathcal{T}_{1-}^{f^{(1)}} = (-\infty, -\beta_2)$, for some $\beta_1, \beta_2 \geq 0$. Also note that when the distortion function is the squared error, the optimal codepoints corresponding to quantization regions $(\mathcal{T}_0^{f^{(1)}}, \mathcal{T}_{1+}^{f^{(1)}}, \mathcal{T}_{1-}^{f^{(1)}})$ are $(\mathbb{E}[X|X \in \mathcal{T}_0^{f^{(1)}}], \mathbb{E}[X|X \in \mathcal{T}_{1+}^{f^{(1)}}], \mathbb{E}[X|X \in \mathcal{T}_{1-}^{f^{(1)}}])$. Hence, we have $D(Q^{(1)}) \leq D(\hat{Q}) \leq D(Q^f)$. □

Based on Lemma 2, we further propose a proposition, and we will prove it later.

**Proposition 1.** *Suppose the source density $p_X$ satisfies Assumption 1. Then, for any communication scheduling policy $f$ satisfying Assumption 2, we can construct a symmetric threshold-based communication scheduling policy $f^{(2)}$ such that*

(1) $\mathcal{T}_0^{f^{(2)}} = (-\beta, \beta), \mathcal{T}_{1+}^{f^{(2)}} = (\beta, \infty), \mathcal{T}_{1-}^{f^{(2)}} = (-\infty, -\beta)$, where $\beta > 0$.

(2) $\mathbb{P}(X \in \mathcal{T}_0^{f^{(2)}}) = \mathbb{P}(X \in \mathcal{T}_0^f)$.

(3) $D(Q^{f^{(2)}}) \leq D(Q^f)$.

To prove Proposition 1, we need to apply results from majorization theory, introduced below. Given a Borel measurable set $A$, we use $A^\sigma$ to denote its symmetric rearrangement, i.e., $A^\sigma = [-a, a]$, and $\mathcal{L}(A^\sigma) = \mathcal{L}(A)$. Given a nonnegative integrable function $p : \mathbb{R} \to \mathbb{R}$, we use $p^\sigma$ to denote its symmetric rearrangement, which is described as follows,

$$p^\sigma(x) := \int_0^\infty \mathbf{1}_{\{z \in \mathbb{R} | p(z) \geq \rho\}^\sigma}(x) d\rho, \quad x \in \mathbb{R}.$$

$\mathbf{1}_{\{z \in \mathbb{R} | p(z) \geq \rho\}^\sigma}(x)$ is the indicator function on whether $x$ is an element of $\{z \in \mathbb{R} | p(z) \geq \rho\}^\sigma$ or not, i.e.,

$$\mathbf{1}_{\{z \in \mathbb{R} | p(z) \geq \rho\}^\sigma}(x) = \begin{cases} 1, & \text{if } x \in \{z \in \mathbb{R} | p(z) \geq \rho\}^\sigma \\ 0, & \text{otherwise} \end{cases}$$

**Definition 1.** *Given two probability densities $p$ and $q$ defined on $\mathbb{R}$, we say $p$ majorizes $q$, denoted by $p \succ q$, if*

$$\int_{|x|<t} q^\sigma(x) dx \leq \int_{|x|<t} p^\sigma(x) dx, \quad \text{for all } t \geq 0$$

**Lemma 3** ([23], Lemma 4). *Let $p_X$ and $p_{X'}$ be probability density functions defined on $\mathbb{R}$. Assume that $p_X$ is even and logconcave, and $p_X \succ p_{X'}$. Then,*

$$\int_{-\infty}^\infty x^2 p_X(x) dx \leq \int_{-\infty}^\infty (x-y)^2 p_{X'}(x) dx, \quad \text{for all } y \in \mathbb{R},$$

*or equivalently,*
$$\text{Var}(X) \leq \text{Var}(X')$$

**Remark 6.** *Lemma 4 in [23] assumes that $p_X$ is even, quasi-concave, and there exists $b \in \mathbb{R}$ such that $p_X$ is non-decreasing on $(-\infty, b]$ and non-increasing on $(b, \infty)$. Note that a positive logconcave function is also quasi-concave. Moreover, it can be easily shown that if $p_X$ is even and logconcave. In view of this, $p_X$ is non-decreasing on $(-\infty, 0]$ and non-increasing on $(0, \infty)$, and hence as stated above Lemma 3 is a valid slightly modified version of Lemma 4 of [23].*

**Lemma 4** ([23], Lemma 2). *Let $p_X$ and $p_{X'}$ be probability density functions defined on $\mathbb{R}$. Assume that $p_X$ is even and logconcave, and $p_X \succ p_{X'}$. Let $A = [-\tau, \tau]$ be any symmetric closed interval such that $\int_A p_X(x)dx > 0$ and let $h : \mathbb{R} \to [0, 1]$ be any function such that $\int_\mathbb{R} h(x)p_{X'}(x)dx = \int_A p_X(x)dx$. Then,*

$$p_{X|X \in A} \succ \frac{h \cdot p_{X'}}{\int_\mathbb{R} h(x) p_{X'}(x) dx},$$

*where*

$$p_{X|X \in A} = \begin{cases} \dfrac{p_X(x)}{\int_A p_X(x)dx}, & \text{if } x \in A \\ 0, & \text{otherwise} \end{cases},$$

*and*

$$\frac{h \cdot p_{X'}}{\int_\mathbb{R} h(x) p_{X'}(x) dx}(x) = \frac{h(x) \cdot p_{X'}(x)}{\int_\mathbb{R} h(x) p_{X'}(x) dx}$$

Combining Lemmas 3 and 4, we have the following lemma.

**Lemma 5.** *Let $p_X$ be an even and logconcave density. Let $A = [-\tau, \tau]$ be any symmetric closed interval such that $\int_A p_X(x)dx > 0$, and let $B$ be any subset of $\mathbb{R}$ such that $\int_B p_X(x)dx = \int_A p_X(x)dx$. Then,*

$$\text{Var}(X|X \in A) \leq \text{Var}(X|X \in B)$$

*Proof.* One can see that $p_X$ majorizes itself. Furthermore, we choose $h(x)$ to be the indicator function on whether $x$ belongs to $B$ or not, i.e., $h(x) = \mathbf{1}_{\{x \in B\}}$. Then, $\int_\mathbb{R} h(x)p_X(x)dx = \int_B p_X(x)dx = \int_A p_X(x)dx$. By Lemma 4, the conditional density $p_{X|X \in A}$ majorizes the conditional density $p_{X|X \in B}$. Since $A$ is symmetric about zero, $p_{X|X \in A}$ is even and logconcave. By Lemma 3, we conclude that $\text{Var}(X|X \in A) \leq \text{Var}(X|X \in B)$. □

To prove Proposition 1, we also need to apply properties of logconcave probability density functions, which is introduced below.

**Lemma 6** ([30], Theorem 6). *Let $p_X$ be a continuously differentiable and logconcave probability density function defined on $(a, b)$. Let $\beta$ be a variable belonging to interval $(a, b)$. Then, the function $G_X(\beta)$, defined below, is monotone decreasing in $\beta$:*

$$G_X(\beta) := \mathbb{E}[X|X > \beta] - \beta. \tag{3}$$

We will frequently refer to this function $G_X(\beta)$ in the rest of the paper. We next provide an extension of Lemma 6 as follows.

**Lemma 7.** *Let $p_X$ be an even and logconcave probability density function defined on $\mathbb{R}$. Furthermore, let $p_X$ be continuously differentiable on $(0, \infty)$, and $\beta$ taking values in $(0, \infty)$. Then, $G_X(\beta)$ as defined by (3) is monotone decreasing in $\beta$.*

*Proof.* Let $Y$ be a random variable such that $Y = |X|$. Denote by $p_Y$ be the probability function of $Y$. Since the probability density of $X$, $p_X$ is even, we have

$$p_Y(y) = \begin{cases} 2p_X(y), & \text{if } y > 0 \\ 0, & \text{otherwise} \end{cases}$$

Since $p_X$ is continuously differentiable on $(0,\infty)$, so is $p_Y$. Furthermore, for any $\beta \in (0,\infty)$, we have

$$\begin{aligned}
\mathbb{E}[Y|Y > \beta] &= \frac{1}{\int_\beta^\infty p_Y(y)dy} \int_\beta^\infty y p_Y(y)dy \\
&= \frac{1}{\int_\beta^\infty 2p_X(y)dy} \int_\beta^\infty y \cdot 2p_X(y)dy \\
&= \frac{1}{\int_\beta^\infty p_X(x)dx} \int_\beta^\infty x p_X(x)dx \\
&= \mathbb{E}[X|X > \beta].
\end{aligned}$$

Hence $G_Y(\beta) = G_X(\beta)$. By Lemma 6, $G_Y(\beta)$ is monotone decreasing in $\beta$. Hence, we conclude that $G_X(\beta)$ is also monotone decreasing in $\beta$ when $\beta > 0$. $\square$

We are now in a position to prove Proposition 1.

*Proof of Proposition 1*: By Lemma 2, given any communication scheduling policy $f^{(0)}$ satisfying Assumption 2, we can construct a threshold-based policy $f^{(1)}$ such that

(1) $\mathcal{T}_0^{f^{(1)}} = (-\beta_2, \beta_1), \mathcal{T}_{1+}^{f^{(1)}} = (\beta_1, \infty), \mathcal{T}_{1-}^{f^{(1)}} = (-\infty, -\beta_2)$.
(2) $\mathbb{P}(X \in \mathcal{T}_i^{f^{(1)}}) = \mathbb{P}(X \in \mathcal{T}_i^{f^{(0)}})$, for all $i \in \{0, 1+, 1-\}$.
(3) $D(Q^{f^{(1)}}) \leq D(Q^{f^{(0)}})$.

Based on policy $f^{(1)}$, we can construct a symmetric threshold-based policy $f^{(2)}$ such that

(1) $\mathcal{T}_0^{f^{(2)}} = (-\beta, \beta), \mathcal{T}_{1+}^{f^{(2)}} = (\beta, \infty), \mathcal{T}_{1-}^{f^{(2)}} = (-\infty, -\beta)$.
(2) $\mathbb{P}(X \in \mathcal{T}_0^{f^{(2)}}) = \mathbb{P}(X \in \mathcal{T}_0^{f^{(1)}})$.

Then, we only need to show that $D(Q^{f^{(2)}}) \leq D(Q^{f^{(1)}})$. Note that $D(Q^{f^{(1)}})$ and $D(Q^{f^{(2)}})$ can be expressed as

$$D(Q^{f^{(1)}}) = \sum_{i \in \{0, 1+, 1-\}} \text{Var}(X|X \in \mathcal{T}_i^{f^{(1)}}) \, \mathbb{P}(X \in \mathcal{T}_i^{f^{(1)}})$$

$$D(Q^{f^{(2)}}) = \sum_{i \in \{0, 1+, 1-\}} \text{Var}(X|X \in \mathcal{T}_i^{f^{(2)}}) \, \mathbb{P}(X \in \mathcal{T}_i^{f^{(2)}})$$

By Lemma 5, $\text{Var}(X|X \in \mathcal{T}_0^{f^{(2)}}) \leq \text{Var}(X|X \in \mathcal{T}_0^{f^{(1)}})$. Since $\mathbb{P}(X \in \mathcal{T}_0^{f^{(2)}}) = \mathbb{P}(X \in \mathcal{T}_0^{f^{(1)}})$, we have

$$\text{Var}(X|X \in \mathcal{T}_0^{f^{(2)}})\mathbb{P}(X \in \mathcal{T}_0^{f^{(2)}}) \leq \text{Var}(X|X \in \mathcal{T}_0^{f^{(1)}})\mathbb{P}(X \in \mathcal{T}_0^{f^{(1)}}).$$

Hence, we will be done if we show that

$$\sum_{i \in \{1+, 1-\}} \text{Var}(X|X \in \mathcal{T}_i^{f^{(2)}}) \, \mathbb{P}(X \in \mathcal{T}_i^{f^{(2)}}) \leq \sum_{i \in \{1+, 1-\}} \text{Var}(X|X \in \mathcal{T}_i^{f^{(1)}}) \, \mathbb{P}(X \in \mathcal{T}_i^{f^{(1)}}).$$

Consider the class of threshold-based communication scheduling policies, denoted by $\mathcal{F}$, whose generic element $f$ is in the form of

$$\mathcal{T}_0^f = (-\gamma_2, \gamma_1), \quad \mathcal{T}_{1+}^f = (\gamma_1, \infty), \quad \mathcal{T}_{1-}^f = (-\infty, -\gamma_2), \quad \gamma_1, \gamma_2 \geq 0,$$

and

$$\mathbb{P}(X \in \mathcal{T}_0^f) = \mathbb{P}(X \in \mathcal{T}_0^{f^{(0)}})$$

It is clear that $f^{(1)}$ and $f^{(2)}$ are elements of $\mathcal{F}$. Let $PD(Q^f)$ be the sum of the mean squared distortions of $Q^f$ over regions $\mathcal{T}_{1+}^f$ and $\mathcal{T}_{1-}^f$, i.e.,

$$\begin{aligned} PD(Q^f) &:= \sum_{i \in \{1+, 1-\}} \mathrm{Var}(X|X \in \mathcal{T}_i^f) \, \mathbb{P}(X \in \mathcal{T}_i^f) \\ &= \mathrm{Var}(X|X < -\gamma_2) \, \mathbb{P}(X < -\gamma_2) + \mathrm{Var}(X|X > \gamma_1) \, \mathbb{P}(X > \gamma_1) \\ &= \mathrm{Var}(X|X > \gamma_2) \, \mathbb{P}(X > \gamma_2) + \mathrm{Var}(X|X > \gamma_1) \, \mathbb{P}(X > \gamma_1), \end{aligned}$$

where the last equality holds since $p_X$ is even. We now show that $f^{(2)}$ is the global minimizer of $PD(Q^f)$ among all elements in $\mathcal{F}$. Since $\mathbb{P}(X \in \mathcal{T}_0^f)$ is a constant, so is $\mathbb{P}(X \in \mathcal{T}_{1+}^f) + \mathbb{P}(X \in \mathcal{T}_{1-}^f)$. Then,

$$\int_{-\infty}^{-\gamma_2} p_X(x)dx + \int_{\gamma_1}^{\infty} p_X(x)dx = \text{constant}$$

Taking the derivatives of both sides with respect to $\gamma_1$, we have

$$\frac{d\gamma_2}{d\gamma_1} \cdot \frac{d}{d\gamma_2} \int_{-a}^{-\gamma_2} p_X(x)dx + \frac{d}{d\gamma_1} \int_{\gamma_1}^{a} p_X(x)dx = 0,$$

which implies that

$$-\frac{d\gamma_2}{d\gamma_1} p_X(-\gamma_2) - p_X(\gamma_1) = 0 \Rightarrow \frac{d\gamma_2}{d\gamma_1} = -\frac{p_X(\gamma_1)}{p_X(\gamma_2)}. \tag{4}$$

The equality above holds because $p_X$ is even. Now taking the derivative of $PD(Q^f)$ with respect to $\gamma_1$, we have

$$\frac{d}{d\gamma_1} PD(Q^f) = \frac{d\gamma_2}{d\gamma_1} \cdot \frac{d}{d\gamma_2} \mathrm{Var}(X|X > \gamma_2) \, \mathbb{P}(X > \gamma_2) + \frac{d}{d\gamma_1} \mathrm{Var}(X|X > \gamma_1) \, \mathbb{P}(X > \gamma_1). \tag{5}$$

The second term in (5) can be computed as follows:

$$\begin{aligned} &\frac{d}{d\gamma_1} \mathrm{Var}(X|X > \gamma_1) \, \mathbb{P}(X > \gamma_1) \\ &= \frac{d}{d\gamma_1} \left( \mathrm{Var}(X|X > \gamma_1) \int_{\gamma_1}^{\infty} p_X(x)dx \right) \\ &= \frac{d}{d\gamma_1} \left( \int_{\gamma_1}^{\infty} x^2 p_X(x)dx - \left(\mathbb{E}[X|X > \gamma_1]\right)^2 \int_{\gamma_1}^{\infty} p_X(x)dx \right) \\ &= \frac{d}{d\gamma_1} \left( \int_{\gamma_1}^{\infty} x^2 p_X(x)dx - \frac{\left(\int_{\gamma_1}^{\infty} x p_X(x)dx\right)^2}{\int_{\gamma_1}^{\infty} p_X(x)dx} \right) \\ &= -\gamma_1^2 p_X(\gamma_1) - \frac{2\int_{\gamma_1}^{\infty} x p_X(x)dx \cdot -\gamma_1 p_X(\gamma_1) \cdot \int_{\gamma_1}^{\infty} p_X(x)dx + \left(\int_{\gamma_1}^{\infty} x p_X(x)dx\right)^2 p_X(\gamma_1)}{\left(\int_{\gamma_1}^{\infty} p_X(x)dx\right)^2} \\ &= -p_X(\gamma_1) \cdot \left(\gamma_1 - \mathbb{E}[X|X > \gamma_1]\right)^2. \end{aligned} \tag{6}$$

Similarly, we have the first term in (5):

$$\frac{d}{d\gamma_2} \mathrm{Var}(X|X > \gamma_2) \, \mathbb{P}(X > \gamma_2) = -p_X(\gamma_2) \cdot \left(\gamma_2 - \mathbb{E}[X|X > \gamma_2]\right)^2 \tag{7}$$

Plugging (4), (6), and (7) into (5), we have

$$\frac{d}{d\gamma_1}PD(Q^f) = -\frac{p_X(\gamma_1)}{p_X(\gamma_2)} \cdot -p_X(\gamma_2) \cdot \big(\gamma_2 - \mathbb{E}[X|X > \gamma_2]\big)^2 - p_X(\gamma_1) \cdot \big(\gamma_1 - \mathbb{E}[X|X > \gamma_1]\big)^2$$

$$= p_X(\gamma_1)\Big[\big(\gamma_2 - \mathbb{E}[X|X > \gamma_2]\big)^2 - \big(\gamma_1 - \mathbb{E}[X|X > \gamma_1]\big)^2\Big] = p_X(\gamma_1)\Big(G_X^2(\gamma_2) - G_X^2(\gamma_1)\Big)$$

By Lemma 7, $G_X(\gamma)$ is a non-negative and monotone decreasing function. Therefore, $G_X^2(\gamma)$ is monotone decreasing, and

$$\frac{d}{d\gamma_1}PD(Q^f) \geq 0, \text{ if } \gamma_1 > \gamma_2, \quad \frac{d}{d\gamma_1}PD(Q^f) = 0, \text{ if } \gamma_1 = \gamma_2, \quad \frac{d}{d\gamma_1}PD(Q^f) \leq 0, \text{ if } \gamma_1 < \gamma_2.$$

which implies that $\gamma_1 = \gamma_2$ (corresponding to $f^{(2)}$) is the global minimizer. $\square$

Now we are in a position to prove Theorem 2 by applying Proposition 1.

*Proof of Theorem 2*: Consider any communication scheduling policy $f$. The expected cost corresponding to $f$ can be computed as follows:

$$J(f) = \mathbb{E}\Big[cU + (X - \hat{X})^2\Big] = \sum_{i \in \{0, 1+, 1-\}} \mathbb{E}\Big[cU + (X - \hat{X})^2 | X \in \mathcal{T}_i^f\Big] \mathbb{P}(X \in \mathcal{T}_i^f)$$

When $X \in \mathcal{T}_0^f$, we have $U = 0$ and $\hat{X} = \mathbb{E}[X | X \in \mathcal{T}_0^f]$. Hence,

$$\mathbb{E}\Big[cU + (X - \hat{X})^2 | X \in \mathcal{T}_0^f\Big] = \mathbb{E}\Big[\big(X - \mathbb{E}[X | X \in \mathcal{T}_0^f]\big)^2 | X \in \mathcal{T}_0^f\Big] = \text{Var}(X | X \in \mathcal{T}_0^f)$$

When $X \in \mathcal{T}_{1+}^f$, we have $U = 1$, $Y = \alpha(1)\big(X - \mathbb{E}[X | X \in \mathcal{T}_{1+}^f]\big)$

$$\hat{X} = \frac{1}{\alpha(1)} \frac{\gamma}{\gamma + 1} \tilde{Y} + \mathbb{E}[X | X \in \mathcal{T}_{1+}^f]$$

$$= \frac{1}{\alpha(1)} \frac{\gamma}{\gamma + 1} (Y + V) + \mathbb{E}[X | X \in \mathcal{T}_{1+}^f]$$

$$= \frac{\gamma}{\gamma + 1} X + \frac{1}{\alpha(1)} \frac{\gamma}{\gamma + 1} V + \frac{1}{\gamma + 1} \mathbb{E}[X | X \in \mathcal{T}_{1+}^f]$$

where

$$\alpha(1) = \sqrt{\frac{P_T}{\text{Var}(X | X \in \mathcal{T}_{1+}^f)}}, \quad \gamma = \frac{P_T}{\sigma_V^2} \quad (8)$$

Hence,

$$\mathbb{E}\big[cU + (X - \hat{X})^2 | X \in \mathcal{T}_{1+}^f\big] = c + \mathbb{E}\bigg[\Big(\frac{1}{\gamma + 1}\big(X - \mathbb{E}[X | X \in \mathcal{T}_{1+}^f]\big) - \frac{1}{\alpha(1)} \frac{\gamma}{\gamma + 1} V\Big)^2 \Big| X \in \mathcal{T}_{1+}^f\bigg]$$

$$= c + \frac{1}{(\gamma + 1)^2} \mathbb{E}\Big[\big(X - \mathbb{E}[X | X \in \mathcal{T}_{1+}^f]\big)^2 | X \in \mathcal{T}_{1+}^f\Big] + \frac{1}{\alpha(1)^2} \frac{\gamma^2}{(\gamma + 1)^2} \mathbb{E}[V^2 | X \in \mathcal{T}_{1+}^f]$$

$$- \frac{2}{\alpha(1)} \frac{\gamma}{(\gamma + 1)^2} \mathbb{E}\Big[\big(X - \mathbb{E}[X | X \in \mathcal{T}_{1+}^f]\big) V | X \in \mathcal{T}_{1+}^f\Big]$$

$$= c + \frac{1}{(\gamma + 1)^2} \text{Var}(X | X \in \mathcal{T}_{1+}^f) + \frac{1}{\alpha(1)^2} \frac{\gamma^2}{(\gamma + 1)^2} \sigma_V^2$$

$$= c + \frac{1}{(\gamma + 1)^2} \text{Var}(X | X \in \mathcal{T}_{1+}^f) + \frac{\gamma}{(\gamma + 1)^2} \text{Var}(X | X \in \mathcal{T}_{1+}^f)$$

$$= c + \frac{1}{\gamma + 1} \text{Var}(X | X \in \mathcal{T}_{1+}^f)$$

The third equality holds since conditioning on the event that $X > 0$, $V$ is zero-mean and independent of $X$ (Assumption 3). The fourth equality holds due to the definitions of $\alpha(1)$ and $\gamma$, which are described by (8). Similarly, one can compute that

$$\mathbb{E}\big[cU + (X - \hat{X})^2 | X \in \mathcal{T}_{1-}^f\big] = c + \frac{1}{\gamma + 1}\text{Var}(X | X \in \mathcal{T}_{1-}^f).$$

Hence, $J(f)$ can be further expressed as

$$J(f) = \text{Var}(X | X \in \mathcal{T}_0^f)\mathbb{P}(X \in \mathcal{T}_0^f) + \frac{1}{\gamma + 1}\text{Var}(X | X \in \mathcal{T}_{1+}^f)\mathbb{P}(X \in \mathcal{T}_{1+}^f) + c\,\mathbb{P}(X \in \mathcal{T}_{1+}^f)$$
$$+ \frac{1}{\gamma + 1}\text{Var}(X | X \in \mathcal{T}_{1-}^f)\mathbb{P}(X \in \mathcal{T}_{1-}^f) + c\,\mathbb{P}(X \in \mathcal{T}_{1-}^f) \quad (9)$$
$$= \frac{1}{\gamma + 1}D(Q^f) + \frac{\gamma}{\gamma + 1}\text{Var}(X | X \in \mathcal{T}_0^f)\mathbb{P}(X \in \mathcal{T}_0^f) + c\,\mathbb{P}(X \in \mathcal{T}_{1+}^f) + c\,\mathbb{P}(X \in \mathcal{T}_{1-}^f)$$

Given any communication scheduling policy $f$, we can construct a symmetric threshold-based communication scheduling policy $f'$ such that

(1) $\mathcal{T}_0^{f'} = (-\beta, \beta), \mathcal{T}_{1+}^{f'} = (\beta, \infty), \mathcal{T}_{1-}^{f'} = (-\infty, -\beta)$.
(2) $\mathbb{P}(X \in \mathcal{T}_0^{f'}) = \mathbb{P}(X \in \mathcal{T}_0^f)$, or equivalently,
$\mathbb{P}(X \in \mathcal{T}_{1+}^{f'}) + \mathbb{P}(X \in \mathcal{T}_{1-}^{f'}) = \mathbb{P}(X \in \mathcal{T}_{1+}^f) + \mathbb{P}(X \in \mathcal{T}_{1-}^f)$

By Proposition 1 and Lemma 5, we have $D(Q^{f'}) \leq D(Q^f)$ and $\text{Var}(X | X \in \mathcal{T}_0^{f'}) \leq \text{Var}(X | X \in \mathcal{T}_0^f)$. Furthermore, we have $\mathbb{P}(X \in \mathcal{T}_0^{f'}) = \mathbb{P}(X \in \mathcal{T}_0^f)$ and $c\,\mathbb{P}(X \in \mathcal{T}_{1+}^{f'}) + c\,\mathbb{P}(X \in \mathcal{T}_{1-}^{f'}) = c\,\mathbb{P}(X \in \mathcal{T}_{1+}^f) + c\,\mathbb{P}(X \in \mathcal{T}_{1-}^f)$. Hence, we conclude that $J(f') \leq J(f)$, which implies that without loss of optimality, we can restrict the search of the optimal communication scheduling policy to the class of symmetric threshold-based policies. Denote by $J(\beta)$ the expected cost corresponding to the symmetric-threshold based communication scheduling policy with threshold $\beta$, where $\beta \geq 0$. By (9), $J(\beta)$ can be computed as

$$J(\beta) = \int_{-\beta}^{\beta} x^2 p_X(x)dx + \frac{1}{\gamma + 1}\text{Var}(X | X < -\beta)\mathbb{P}(X < -\beta) + c\int_{-\infty}^{-\beta} p_X(x)dx$$
$$+ \frac{1}{\gamma + 1}\text{Var}(X | X > \beta)\mathbb{P}(X > \beta) + c\int_{\beta}^{\infty} p_X(x)dx$$
$$= 2\int_0^{\beta} x^2 p_X(x)dx + 2\frac{1}{\gamma + 1}\text{Var}(X | X > \beta)\mathbb{P}(X > \beta) + 2c\int_{\beta}^{\infty} p_X(x)dx,$$

where the second equality holds since $p_X$ is even. Taking the derivative of $J(\beta)$ with respect to $\beta$, and by (6), we have

$$\frac{d}{d\beta}J(\beta) = 2p_X(\beta)\Big(\beta^2 - \frac{1}{\gamma + 1}\big(\mathbb{E}[X | X > \beta] - \beta\big)^2 - c\Big)$$
$$= 2p_X(\beta)\Big(\beta^2 - \frac{1}{\gamma + 1}G_X^2(\beta) - c\Big)$$

Since $c > 0$ and $G_X(\beta)$ is monotone decreasing, there exists a unique $\beta^*$ in $[0, \infty)$ such that

$$\beta^{*2} = \frac{1}{\gamma + 1}G_X^2(\beta^*) + c$$

Furthermore, $dJ(\beta)/d\beta < 0$ when $\beta < \beta^*$ and $dJ(\beta)/d\beta > 0$ when $\beta > \beta^*$. Hence, $\beta^*$ is the unique global minimizer among all $\beta > 0$. □

**Remark 7.** *If the density function $p_X$ has support $(-a, a)$ and $0 < a < \beta^*$, then $dJ(\beta)/d\beta$ is always negative, which implies that the minimizing $\beta$ is just $a$. This means the optimal communication scheduling policy is to always choose not to transmit regardless of sensor's observation. Such a case can occur when the communication cost is very high.*

## 4 Optimization Problem with hard constraints

To present our main results for the problem with hard constraints, we introduce a number of terms. First, we let $E_t$ denote the number of communication opportunities left at the beginning of the $t$-th time interval, i.e., $E_t = N - \sum_{i=1}^{t-1} U_i$. Then, evolution of $E_t$ is described by

$$
\begin{aligned}
E_1 &= N \\
E_t &= E_{t-1} - U_{t-1}, \quad t \geq 2
\end{aligned}
\tag{10}
$$

Furthermore, the communication constraint is

$$U_t \leq E_t, \quad \text{for all } t = 1, 2, \ldots, T \tag{11}$$

Recall that $U_{1:t-1}$ is the common information shared by the sensor, the encoder, and the decoder, and hence $E_t$ is also known by all the decision makers.

Second, we let $J^*(t, E_t)$ be the optimal cost to go if the system is initialized at time $t$ (or equivalently, at the beginning of the $t$-th time interval) with $E_t$ number of communication opportunities. Specifically, we define $J^*(T+1, \cdot) = 0$ for any number of communication opportunities.

Third, for any $E_t > 0$, we let $c(t, E_t)$ denote the difference between two optimal costs to go, i.e.,

$$c(t, E_t) = J^*(t+1, E_t - 1) - J^*(t+1, E_t).$$

**Remark 8.** *$c(t, E_t)$ can be interpreted as the opportunity cost for choosing to communicate with the estimator rather than not to communicate.*

The following theorem ensures that without loss of optimality we can restrict all the decision makers to consider only their current inputs and $E_t$ when making decisions at time $t$. Furthermore, the optimal decision policies can be obtained via solving the dynamic programming equation.

**Theorem 3.** *Consider the optimization problem with hard constraint as formulated in section 2-D. Without loss of optimality, we can restrict the communication scheduling, encoding and decoding policies to the forms:*

$$U_t = f_t(X_t, E_t), \ Y_t = g_t(\tilde{X}_t, E_t), \ \hat{X}_t = h_t(\tilde{Y}_t, S_t, E_t)$$

*Furthermore, the optimal decision policies $(f_t^*, g_t^*, h_t^*)$ together with the optimal cost to go $J^*(t, E_t)$ can be obtained by solving the dynamic programming (DP) equation:*

$$
\begin{aligned}
J^*(T+1, \cdot) &= 0 \\
J^*(t, E_t) &= \min_{f_t, g_t, h_t} \mathbb{E}\left\{(X_t - \hat{X}_t)^2 + J^*(t+1, E_{t+1})\right\}
\end{aligned}
\tag{12}
$$

The proof of Theorem 3 is similar to that of Theorem 1, and hence we keep using notations $(\mathbf{f}_{a:b}, \mathbf{g}_{a:b}, \mathbf{h}_{a:b})$ and $(I_{st}, I_{et}, I_{dt})$, and $I_t$.

*Proof.* It is easy to see validity of the following equality:

$$
\begin{aligned}
&\min_{\mathbf{f}, \mathbf{g}, \mathbf{h}} J(\mathbf{f}, \mathbf{g}, \mathbf{h}) \\
&= \min_{f_1, g_1, h_1} \mathbb{E}\bigg\{cU_1 + (X_1 - \hat{X}_1)^2 + \min_{f_2, g_2, h_2} \mathbb{E}\bigg\{cU_2 + (X_2 - \hat{X}_2)^2 + \ldots \\
&\quad + \min_{f_{T-1}, g_{T-1}, h_{T-1}} \mathbb{E}\bigg\{cU_{T-1} + (X_{T-1} - \hat{X}_{T-1})^2 + \min_{f_T, g_T, h_T} \mathbb{E}\big\{cU_T + (X_T - \hat{X}_T)^2\big\}\bigg\}\ldots\bigg\}
\end{aligned}
\tag{13}
$$

Then, at time $t = T$, the optimization problem, call it *Problem T1*, is to design $(f_T, g_T, h_T)$ minimizing

$$J_T(f_T, g_T, h_T) := \mathbb{E}\left\{(X_T - \hat{X}_T)^2\right\}$$

where

$$U_T = f_T(X_T, I_{sT}), Y_T = g_T(\tilde{X}_T, I_{eT}), \hat{X}_T = h_T(\tilde{Y}_T, S_T, I_{dT}).$$

Consider another problem, call it *Problem T2*, where $I_T$ is available to the sensor, the encoder, and the decoder, and one needs to design $(f'_T, g'_T, h'_T)$ minimizing

$$J_{T_2}(f'_T, g'_T, h'_T) := \mathbb{E}\left\{cU_T + (X_T - \hat{X}_T)^2\right\},$$

where

$$U_T = f'_T(X_T, I_T), \quad Y_T = g'_T(\tilde{X}_T, I_T), \quad \hat{X}_T = h'_T(\tilde{Y}_T, S_T, I_T).$$

Since the sensor, the encoder, and the decoder can always ignore the redundant information and behave as if they only know $I_{sT}, I_{eT}, I_{dT}$, respectively, the performance of the system in *Problem T2* is no worse than that in *Problem T1*, i.e.,

$$\min_{(f'_T, g'_T, h'_T)} J_{T_2}(f'_T, g'_T, h'_T) \leq \min_{(f_T, g_T, h_T)} J_{T_1}(f_T, g_T, h_T).$$

Similarly, consider a third problem, call it *Problem T3*, where $I_{sT}, I_{eT}, I_{dT}$ are not available to the sensor, the encoder, and the decoder, respectively. Instead, all decision makers share the information about $E_T$. One needs to design $(f''_T, g''_T, h''_T)$ minimizing

$$J_{T_3}(f''_T, g''_T, h''_T) = \mathbb{E}\left\{cU_T + (X_T - \hat{X}_T)^2\right\},$$

where

$$U_T = f''_T(X_T, E_T), \quad Y_T = g''_T(\tilde{X}_T, E_T), \quad \hat{X}_T = h''_T(\tilde{Y}_T, S_T, E_T)$$

Since the sensor, the encoder, and the decoder can deduce $E_T$ from $I_{sT}, I_{eT}, I_{dT}$, respectively, the system in *Problem T1* cannot perform worse than the system in *Problem T3*. Hence,

$$\min_{(f_T, g_T, h_T)} J_{T_1}(f_T, g_T, h_T) \leq \min_{(f''_T, g''_T, h''_T)} J_{T_3}(f''_T, g''_T, h''_T).$$

We now consider system 2. Note that the random variables $(X_T, V_T)$, the distortion function $\rho(\cdot, \cdot)$, and the power constraint of the encoder do not depend on $I_T$. Furthermore, the communication constraint as described by (11) depends on $I_T$ only via $E_T$. Therefore, there is no loss of optimality in restricting

$$U_T = f'_T(X_T, E_T), Y_T = g'_T(\tilde{X}_T, E_T), \hat{X}_T = h'_T(\tilde{Y}_T, S_T, E_T),$$

and thus

$$\min_{(f'_T, g'_T, h'_T)} J_{T_2}(f'_T, g'_T, h'_T) = \min_{(f''_T, g''_T, h''_T)} J_{T_3}(f''_T, g''_T, h''_T)$$

The equality above indicates that in *Problem T1*, the sensor, the encoder, and the decoder can ignore their information about the past, namely $I_{sT}, I_{eT}$, and $I_{dT}$, respectively, and there is no loss of optimality in restricting

$$U_T = f_T(X_T, E_T), Y_T = g_T(\tilde{X}_T, E_T), \hat{X}_T = h_T(\tilde{Y}_T, S_T, E_T).$$

Equivalently, all the decision makers in system 1 can behave as if the system is initialized at time $T$ with $E_T$ number of communication opportunities. Hence,

$$\min_{f_T, g_T, h_T} J_T(f_T, g_T, h_T) = J^*(T, E_T).$$

Now consider the optimization problem at stage $T - 1$. By (13), the problem is to design $(f_{T-1}, g_{T-1}, h_{T-1})$ minimizing

$$J_{T-1}(f_{T-1}, g_{T-1}, h_{T-1}) := \mathbb{E}\left\{(X_{T-1} - \hat{X}_{T-1})^2 + \min_{f_T, g_T, h_T} \mathbb{E}\left\{\sum_{t=1}^T cU_t + (X_t - \hat{X}_t)^2\right\}\right\}$$

$$= \mathbb{E}\left\{(X_{T-1} - \hat{X}_{T-1})^2 + J^*(T, E_T)\right\} \quad (14)$$

$$= \mathbb{E}\left\{(X_{T-1} - \hat{X}_{T-1})^2 + J^*(T, E_{T-1} - U_{T-1})\right\}$$

Depending on whether $E_{T-1}$ is zero or not, we have the following: when $E_{T-1} = 0$, by the communication constraint, $U_T = U_{T-1} = 0$. Then, $J^*(T, 0)$ is a constant, which is independent of $(f_{T-1}, g_{T-1}, h_{T-1})$. Hence, minimizing $J_{T-1}(f_{T-1}, g_{T-1}, h_{T-1})$ is equivalent to minimizing $\mathbb{E}\{(X_{T-1} - \hat{X}_{T-1})^2\}$. When $E_{T-1} > 0$, (14) can be further expressed as

$$J_{T-1}(f_{T-1}, g_{T-1}, h_{T-1}) = J^*(T, E_{T-1}) + \mathbb{E}\left\{(X_{T-1} - \hat{X}_{T-1})^2 + c(T, E_{T-1}) \cdot U_{T-1}\right\}$$

Note that $J^*(T, E_{T-1})$ is independent of $(f_{T-1}, g_{T-1}, h_{T-1})$. Hence, minimizing $J_{T-1}(f_{T-1}, g_{T-1}, h_{T-1})$ is equivalent to minimizing $\mathbb{E}\{(X_{T-1} - \hat{X}_{T-1})^2 + c(T, E_{T-1})U_{T-1}\}$. In both cases, the random variables $(X_{T-1}, V_{T-1})$, the distortion function $\rho(\cdot, \cdot)$, and the power constraint of the encoder do not depend on $I_{T-1}$. Furthermore, (in both cases) the communication constraints and (in the second case) the opportunity cost $c(T - 1, E_{T-1})$ depend on $I_{T-1}$ only via $E_{T-1}$. By a similar argument as above, it is sufficient for all the decision makers to consider only $E_{T-1}$ instead of $I_{s(T-1)}, I_{e(T-1)}, I_{d(T-1)}$. Equivalently, all the decision makers can behave as if the system is initialized at time $T - 1$ with $E_{T-1}$ number of communication opportunities. Hence,

$$\min_{f_{T-1}, g_{T-1}, h_{T-1}} J_{T-1}(f_{T-1}, g_{T-1}, h_{T-1}) = J^*(T - 1, E_{T-1}).$$

By induction one can show that there is no loss of generality by restricting $U_t = f_t(X_t, E_t), Y_t = g_t(\tilde{X}_t, E_t), \hat{X}_t = h_t(\tilde{Y}_t, S_t, E_t)$, and the optimal cost to go satisfies the DP equation (12). $\square$

Let us take a closer look at the DP equation (12). When $E_t = 0$, by the communication constraint $U_t = f_t(X_t, E_t) = 0$ regardless of the realization of $X_t$. Consequently, we have $E_{t+1} = 0$. Then, the DP equation can be easily solved as follows:

$$J^*(t, 0) = \min_{f_t, g_t, h_t} \mathbb{E}\left\{(X_t - \hat{X}_t)^2\right\} + J^*(t + 1, 0) = \mathrm{Var}(X_t) + J^*(t + 1, 0)$$

The last equality holds since without any information about $X_t$, the optimal estimator is $\mathbb{E}[X_t]$ and the mean squared error is $\mathrm{Var}(X_t)$. When $E_t > 0$, the DP equation can be written as

$$J^*(t, E_t) = \min_{f_t, g_t, h_t} \mathbb{E}\left\{(X_t - \hat{X}_t)^2 + J^*(t + 1, E_{t+1})\right\} = J^*(t, E_t) + \min_{f_t, g_t, h_t} \mathbb{E}\left\{c(t, E_t)U_t + (X_t - \hat{X}_t)^2\right\} \quad (15)$$

Note that the minimization in the second line of (15) is just the one-stage problem discussed in section 3 with communication cost $c(t, E_t)$. This now motivates us to make the following assumptions.

**Assumption 5.** *The sensor is restricted to apply the communication scheduling policies $f_t$ such that for any $1 \leq t \leq T$ and $E_t > 0$,*

$$\mathbb{E}[X_t | f_t(X_t, E_t) = 1, X_t < 0] < \mathbb{E}[X_t | f_t(X_t, E_t) = 0] < \mathbb{E}[X_t | f_t(X_t, E_t) = 1, X_t > 0]$$

**Assumption 6.** *The encoder and the decoder are restricted to apply piecewise affine policies:*

$$g_t(\tilde{X}_t, E_t) = \begin{cases} S_t \cdot \alpha(S_t) \cdot (X_t - \mathbb{E}[X_t | U_t = 1, S_t]), & \text{if } U_t = 1 \\ 0, & \text{if } U_t = 0 \end{cases}$$

$$h_t(\tilde{Y}_t, S_t, E_t) = \begin{cases} S_t \cdot \dfrac{1}{\alpha(S_t)} \dfrac{\gamma}{\gamma + 1} \tilde{Y}_t + \mathbb{E}[X_t | U_t = 1, S_t], & \text{if } U_t = 1 \\ \mathbb{E}[X_t | U_t = 0], & \text{if } U_t = 0 \end{cases}$$

where $\gamma = P_T/\sigma_V^2$, $\alpha(S_t) = \sqrt{P_T/\text{Var}(X_t|U_t = 1, S_t)}$, and $\text{Var}(X_t|U_t = 1, S_t)$ is the conditional variance.

Then, we have the following theorem on the optimality of symmetric threshold-based communication scheduling strategy.

**Theorem 4.** *Consider the problem with hard constraints under Assumptions 1, 3, 5 and 6*

$$f_t^*(X_t, E_t) = \begin{cases} 1, & \text{if } E_t > 0 \text{ and } |X_t| > \beta_t^*(E_t) \\ 0, & \text{if } E_t = 0 \text{ or } |X_t| \leq \beta_t^*(E_t) \end{cases} \quad (16)$$

*where $\beta_t^*(E_t)$ is non-negative and is the unique solution to the fixed-point equation:*

$$\beta^2 = \frac{1}{\gamma+1} G_X^2(\beta) + c(t, E_t), \quad G_X(\beta) := \mathbb{E}[X_t | X_t > \beta] - \beta, \quad \beta \geq 0 \quad (17)$$

**Remark 9.** *The major differences between the problem considered in this paper and the problems considered in [31] and [32] are as follows: In [31] and [32], we had restricted the sensor to apply symmetric threshold-based policies and shown that the optimal encoding and decoding policies are piecewise affine. Furthermore, the results only hold for specific source and noise densities (e.g., Laplace source and Gamma noise with specific parameters). In this paper, however, we restrict the encoder and the decoder to apply piecewise affine encoding and decoding policies, and we show that under some weak technical assumptions (Assumptions 2 and 5), the optimal communication scheduling policy is symmetric and threshold-based. Moreover, the results hold for a large class of source densities (e.g., general even and logconcave densities).*

**Remark 10.** *Consider the case where $E_t > T - t$, that is, the sensor is always allowed to communicate with the estimator for the remaining time. First, we note that the opportunity cost $c(t, E_t)$ is zero. Since $G_X(0) = \mathbb{E}[X | X > 0] > 0$, the solution to (17) is not zero, which means it is strictly positive. Then, even though the sensor can always communicate with the estimator, the optimal communication policy is still the threshold based policy with threshold $\beta_t^*(E_t) > 0$, which might seem counter-intuitive: why would the sensor not transmit its observation although it is allowed to do so? This surprising result is due to the fact that threshold information, i.e., whether or not the state sample belongs to a fixed, known interval, might be more informative than a noisy observation of the state at the output of the noisy channel. Hence, it might be better not to communicate explicitly over the noisy channel but rely on the side channel which signals where the sample lies. For example, at the extreme case of a very noisy channel ($\gamma \to 0$) the output of the communication channel, $\tilde{Y}_t$, is effectively useless, irrespective of the realization $X_t$. However, depending on the threshold and the realization $X_t$, thresholding information could be significantly more informative.*

## 5 Numerical Results

In this section, we present the numerical results for the problem with hard constraint. We select the source density to be Laplace density with parameter $\lambda$, i.e.,

$$p_X(x) = \begin{cases} \frac{1}{2}\lambda e^{-\lambda x}, & \text{if } x \geq 0 \\ \frac{1}{2}\lambda e^{\lambda x}, & \text{if } x < 0 \end{cases}$$

Then, it is easy to see that

$$G_X(\beta) = \mathbb{E}[X_t | X_t > \beta] - \beta = \frac{1}{\lambda}, \quad \text{for all } \beta \geq 0.$$

Hence, the solution to (17) is

$$\beta_t^*(E_t) = \sqrt{\frac{1}{\gamma+1}\frac{1}{\lambda^2} + c(t, E_t)} = \sqrt{m + c(t, E_t)},$$

where

$$m := \frac{1}{\gamma+1}\frac{1}{\lambda^2}.$$

Then, the optimal communication scheduling policy is described by (16). Furthermore, the optimal encoding/decoding policies $(g_t^*, h_t^*)$ are as follows

$$g_t(\tilde{X}_t, E_t) = \begin{cases} \alpha \cdot \big(|X_t| - \beta_t^*(E_t) - \lambda^{-1}\big), & \text{if } U_t = 1 \\ 0, & \text{if } U_t = 0 \end{cases}$$

$$h_t(\tilde{Y}_t, S_t, E_t) = \begin{cases} S_t \cdot \Big(\dfrac{1}{\alpha}\dfrac{\gamma}{\gamma+1}\tilde{Y}_t + \beta_t^*(E_t) + \lambda^{-1}\Big), & \text{if } U_t = 1 \\ 0, & \text{if } U_t = 0 \end{cases}$$

where $\gamma = P_T/\sigma_V^2$, and $\alpha = \sqrt{P_T/\lambda^{-2}}$. By plugging the optimal communication scheduling, encoding, and decoding policies $(f_t^*, g_t^*, h_t^*)$ into the DP equation (12), we obtain the explicit update rule for the optimal cost to go $J^*(t, E_t)$, as shown below

$$\begin{aligned}
J^*(t, E_t) &= J^*(t+1, E_t) + 2\lambda^{-2}, & \text{if } E_t = 0 \\
J^*(t, E_t) &= J^*(t+1, E_t) + 2\lambda^{-2} - 2\big(\beta_t^*(E_t)\lambda^{-1} + \lambda^{-2}\big)e^{-\lambda\beta_t^*(E_t)}, & \text{if } E_t > 0
\end{aligned} \quad (18)$$

We choose the parameters as follows: $T = 100$, $\lambda = 1$, and the signal-to-noise ratio (SNR) $\gamma = 0.1, 1, 10$. We solve the optimal cost to go $J^*(t, E_t)$ by applying the update rule (18). We have plotted the optimal 100-stage estimation error $J^*(1, N)$ versus the number of communication opportunities $N$ under different SNRs, as shown in Fig. 2.

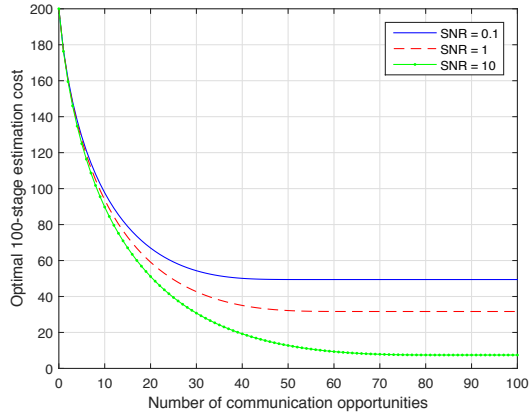

Fig. 2. 100-stages estimation error vs. number of communication opportunities

One can see that, as to be expected, for each fixed SNR, the optimal 100-stage estimation error is non-increasing in terms of the number of communication opportunities. To be more specific, there exists a threshold on the number of communication opportunities (call it *opportunity threshold*) such that the optimal 100-stage estimation error decreases when the number of communication opportunities is below the threshold, and it stays constant above the threshold. We call *minimal error* as the optimal 100-stage estimation error with the number of communication opportunities above the opportunity threshold. One can also see from Fig. 2 that when the SNR increases, the opportunity threshold increases and the minimal error decreases. The existence of opportunity threshold can be interpreted as follows: since we restrict the sensor to apply the threshold based policy with threshold $\beta_t^*(E_t) = \sqrt{c(t, E_t) + m} \geq \sqrt{m}$, the expected number of communication opportunities that will be used is upper bounded by $T \cdot \mathbb{P}(|X_t| \geq \sqrt{m}) = Te^{-\lambda\sqrt{m}}$. Therefore, when the communication opportunities is greater than $Te^{-\lambda\sqrt{m}}$, the optimal expected estimation error will not decrease even though the sensor can have more communication opportunities. It can also be checked from Fig. 2 that the opportunity thresholds under different signal to noise ratios are roughly $Te^{-\lambda\sqrt{m}}$. Moreover, since $m = \frac{1}{\gamma+1}\frac{1}{\lambda^2}$, $Te^{-\lambda\sqrt{m}} = Te^{-1/\sqrt{\gamma+1}}$, which is an increasing function of the SNR $\gamma$. Hence, the opportunity threshold increases with the SNR.

Fig. 3 depicts a sample path of the number of communication opportunities left when the sensor applies the symmetric threshold-based scheduling policies described in Theorem 4. When generating the plot, we have chosen $T = 100$, $\lambda = 1$, $\gamma = 0.1$, and the number of communication opportunities $N = 50$. One can see that the communication opportunities are not used up by the end of the time horizon.

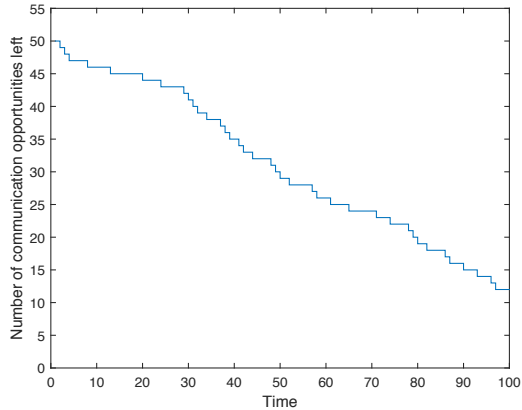

Fig. 3. A sample path of the number of communication opportunities left vs. time

When the number of communication opportunities is larger than the opportunity threshold, the optimal estimation error does not change with respect to the number of communication opportunities. Without loss of generality, we can assume that the sensor is allowed to communicate at each step, that is, $N = T$. Then, the opportunity cost is $c(t, E_t) = 0$. Recall that $\beta_t^*(E_t) = \sqrt{c(t, E_t) + m}$ and $m = \frac{1}{\gamma+1}\frac{1}{\lambda^2}$. Hence, the update rule for the cost function can be simplified as follows:

$$J^*(t, T) = J^*(t+1, T) + \left(\frac{2}{\lambda^2} - \left(\frac{2\sqrt{m}}{\lambda} + \frac{2}{\lambda^2}\right) \cdot e^{-\lambda\sqrt{m}}\right)$$

with $J^*(T+1, T) = 0$, which implies that

$$J^*(1, T) = T \cdot \left(\frac{2}{\lambda^2} - \left(\frac{2\sqrt{m}}{\lambda} + \frac{2}{\lambda^2}\right) \cdot e^{-\lambda\sqrt{m}}\right) = T \cdot 2\lambda^{-2} \cdot \left[1 - \left(\frac{1}{\sqrt{1+\gamma}} + 1\right) \cdot e^{-\frac{1}{\sqrt{1+\gamma}}}\right].$$

It is straightforward to check that $J^*(1, T)$ is a decreasing function of the SNR $\gamma$. Hence, the minimal error decreases as the SNR increases.

Plotting the opportunity threshold $Te^{-\lambda m}$ versus minimal error $J^*(1, T)$ under different SNRs (dash line) in Fig. 2, we arrive at Figure 4. One can see that the intersection between the dash line and each solid line is roughly the turning point of the solid line. Therefore, the plot of opportunity threshold versus minimal error under different SNRs is an important one. In fact, the plot suggests the lowest capacity of the battery that one should choose when building a physical system so that the expected estimation error is minimized. In addition, the plot predicts the minimal expected estimation error.

Consider the asymptotic case where the SNR $\gamma \to \infty$, and thus $m = \frac{1}{\gamma+1}\frac{1}{\lambda^2} \to 0$. Then the opportunity threshold $Te^{-\lambda m} \to T$, and the minimal error $J^*(1, T) \to 0$. Hence, the optimal 100-stage estimation error will be strictly decreasing in terms of the number of communication opportunities in the asymptotic case, as also noted in the prior work [12]. Moreover, the estimation error will reach zero when the number of communication opportunities is equal to the time horizon.

## 6 Conclusions

In this paper, we have considered a communication scheduling and remote estimation problem with a noisy communication channel. Under some technical conditions, we have obtained optimal solutions for both soft-constrained and hard-constrained problems, which consist of a symmetric threshold-based communication scheduling strategy and a pair of piecewise affine encoding/decoding strategies. Moreover, we have generated numerical results to illustrate the effect of the presence of channel noise.

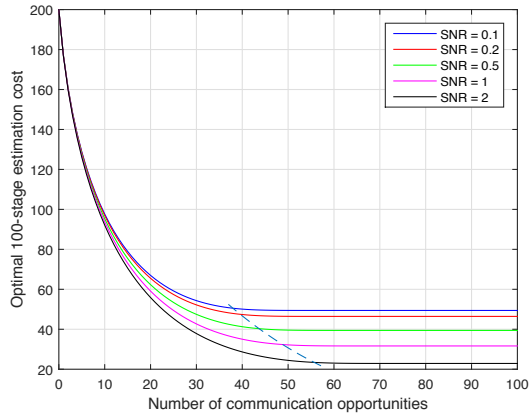

Fig. 4. Opportunity threshold vs. minimal error under different signal to noise ratios

There are several directions for future work: (1) Here, we assumed that the sensor is restricted to apply the communication scheduling policy satisfying Assumption 2. Under this assumption, we showed that without loss of optimality, one can restrict the search of optimal communication scheduling policy to the symmetric class. However, it is plausible that this assumption can be relaxed or even removed, that is, for any communication scheduling policy (which may or may not satisfy Assumption 2), there exists a symmetric policy achieving no greater costs. This possible extension is not immediate and requires more effort. (2) Here, we considered the setting with a noisy channel. It will be interesting to consider a more general setting where there are two channels. One is noisy but not costly, and the other one is perfect (has high communication quality) but is costly. Then, at each time step, the sensor needs to choose whether to transmit its observation or not. If the sensor chooses to transmit its oberservation, it also needs to choose which channel it will use. More details on this problem can be found in [33]. (3) Here, the encoding power was taken to be time invariant. What if the encoder can distribute its total encoding power over the time horizon? More details on problems with power allocation can be found in [34, 35]. (4) Finally, we considered here a one-dimensional system. It would be interesting to consider extensions to multi-dimensional systems, that is, when the input and the communication channel are chosen from multi-dimensional spaces.